
\newcount\capitolo
\newcount\incremento
\newcount\paragrafo

\font\teneuf=eufm10
\font\seveneuf=eufm7
\font\fiveeuf=eufm5
\newfam\euffam
\textfont\euffam=\teneuf
\scriptfont\euffam=\seveneuf
\scriptscriptfont\euffam=\fiveeuf
\def\got{\fam\euffam}

\font\tenmsy=msym10
\font\sevenmsy=msym7
\font\fivemsy=msym5
\font\tenmsx=msxm10
\font\sevenmsx=msxm7
\font\fivemsx=msxm5
\newfam\msxfam
\newfam\msyfam
\def\bbb{\fam\msyfam}
\def\Bbb#1{{\bbb{#1}}}
\textfont\msyfam=\tenmsy
\scriptfont\msyfam=\sevenmsy
\scriptscriptfont\msyfam=\fivemsy
\textfont\msxfam=\tenmsx
\scriptfont\msxfam=\sevenmsx
\scriptscriptfont\msxfam=\fivemsx
\def\hexnumber@#1{\ifcase#1 0\or1\or2\or3\or4\or5\or6\or7\or8\or9\or
	A\or B\or C\or D\or E\or F\fi }
\edef\msx@{\hexnumber@\msxfam}
\edef\msy@{\hexnumber@\msyfam}

\font\eightrm=cmr8
\font\eightbf=cmbx8
\font\eightit=cmti8
\font\eightsl=cmsl8
\font\eightmus=cmmi8
\def\smalltype{\let\rm=\eightrm \let\bf=\eightbf \let\it=\eightit
\let\sl=\eightsl \let\mus=\eightmus \baselineskip=9.5pt minus .75pt
\rm}

\def\hexnumber@#1{\ifcase#1 0\or1\or2\or3\or4\or5\or6\or7\or8\or9\or
	A\or B\or C\or D\or E\or F\fi }
\edef\msy@{\hexnumber@\msyfam}
\mathchardef\ltimes="2\msy@6E
\mathchardef\twoheadrightarrow="3\msx@10

 3

\def\para#1{\global\advance\paragrafo1%
\par\bigbreak\centerline{\bf\number\paragrafo.\ \bf#1}%
\global\incremento=0%
\par\medbreak}
\def\df{\global\advance\incremento by1 {\bf Definition.}}
\def\pf{{\it Proof.\ \ }}
\def\prop{\global\advance\incremento by1 {\bf Proposition.
}}
\def\th{\global\advance\incremento by1 {\bf Theorem.
}}
\def\rem{\global\advance\incremento by1 {\bf Remark.
}}
\def\rems{\global\advance\incremento by1 {\bf Remarks.
}}
\def\lem{\global\advance\incremento by1 {\bf Lemma.
}}
\def\cor{\global\advance\incremento by1 {\bf Corollary.
}}

\def\pr{^{\prime}}
\def\pp{{\prime}}
\def\ppr{^{\prime\prime}}

\def\A{{\rm A}}

\def\R{{\rm R}}

\def\Q{{\rm Q}}

\def\Y{{\rm Y}}
\def\X{{\rm X}}
\def\T{{\rm T}}

\def\l{{\lambda}}
\def\o{{\omega}}

\def\d{{\delta}}
\def\si{{\sigma}}
\def\t{{\theta}}

\def\vir{,\ldots ,}
\def\dd{\displaystyle}

\def\P {{\rm P}}

\def\cvd{\hfill $\sqcap \hskip-6.5pt \sqcup$}        
\overfullrule=0pt                                    

\footline={\hss\lower .8cm \hbox{\tenrm\folio}\hss}
\hsize=15 true cm
\vsize=20.5 true cm
\baselineskip=.5cm
\parindent=0pt
\topskip=0cm
\voffset 1.5 truecm


\def\cb{{\bf c}}
\def\db{{\bf d}}
\def\eb{{\bf e}}
\def\fb{{\bf f}}

\def\kb{{\bf k}}

\def\qb{{\bf q}}

\def\tb{{\bf t}}
\def\vb{{\bf v}}
\def\wb{{\bf w}}
\def\xb{{\bf x}}
\def\yb{{\bf y}}
\def\Hb{{\bf H}}
\def\Ub{{\bf U}}
\def\Yb{{\bf Y}}

\def\ot{{\otimes}}
\def\K{{\kb}}

\def\E{{\eb}}
\def\F{{\fb}}
\def\X{{\rm X}}
\def\Y{{\rm Y}}
\def\M{{\rm M}}

\def\V{{\bf V}}
\def\interi{{\ZZ}}

\def\n{{\NN}}

\def\CC{{\Bbb C}}
\def\NN{{\Bbb N}}
\def\ZZ{{\Bbb Z}}
\def\Ac{{\cal A}}
\def\Bc{{\cal B}}
\def\Mc{{\cal M}}
\def\Nc{{\cal N}}
\def\ti{{\T_i}}
\def\tii{{\T_i^{^{-1}}}}
\def\Htor{{\ddot\Hb_{_\Ac}}}
\def\Haf{{\dot\Hb_{_\Ac}}}
\def\Hafa{{\dot\Hb_{_\Ac}^{^{(1)}}}}
\def\Hafb{{\dot\Hb_{_\Ac}^{^{(2)}}}}
\def\Hf{{\Hb_{_\Ac}}}
\def\htor{{\ddot\Hb}}
\def\haf{{\dot\Hb}}
\def\hafa{{\dot\Hb^{^{(1)}}}}
\def\hafb{{\dot\Hb^{^{(2)}}}}
\def\hf{{\Hb}}
\def\Uf{{\bf U}_{_\Bc}}
\def\Uaf{\dot{\bf U}_{_\Bc}}
\def\Uafa{\dot{\bf U}_{_\Bc}^{^{(1)}}}
\def\Uafb{\dot{\bf U}_{_\Bc}^{^{(2)}}}
\def\Utor{\ddot{\bf U}_{_\Bc}}
\def\utor{\ddot{\bf U}}
\def\uaf{\dot{\bf U}}
\def\uafa{\dot{\bf U}^{^{(1)}}}
\def\uafb{\dot{\bf U}^{^{(2)}}}
\def\uf{{\bf U}}
\def\vu{{\bf v}}
\def\kt{\kb_{_\t}}
\def\et{{\eb_{_\t}}}
\def\ft{{\fb_{_\t}}}
\def\modulo{M\ot_{_\hf}\V^{^{\otimes l}}}
\def\jb{{\bf j}}
\def\jbp{{\bf j_1}}
\def\jbpp{{\bf j_2}}
\def\jbppp{{\bf j_3}}
\def\tbp{{\bf j_4}}

\def\set#1{{[#1]^{^{\times}}}}

\centerline{{\bf  SCHUR DUALITY IN THE TOROIDAL SETTING}}

\vskip 12mm
\centerline {
{\bf M.Varagnolo}
\footnote{$^{^1}$}{\smalltype Dipartimento di Matematica,
via della Ricerca Scientifica, 00133 Roma, Italy.
varagnolo@vax.mat.utovrm.it},
{\bf E.Vasserot}
\footnote{$^{^2}$}{\smalltype Ecole Normale Sup\'erieure, 45, rue d'Ulm,
75005 Paris, France.
vasserot@dmi.ens.fr}
}

\vskip 12mm

\hskip1cm The classical Frobenius-Schur duality gives a correspondence
between finite dimensional representations
of the symmetric and the linear groups. The goal of the present paper
is to extend this construction to the quantum toroidal setup with only
elementary (algebraic) methods. This work can be seen as a
continuation of [J], [D1] and [C2] (see also [CP] and [GRV])
where the cases of the quantum groups $\Ub_q\,(\got{sl}\,(n))$,
$\Yb\,({\widehat{\got{sl}}}\,(n))$ (the Yangian) and
$\Ub_q\,({\widehat{\got{sl}}}\,(n))$ are given.
In the toroidal setting the two algebras involved are a deformation
of Cherednik's double affine Hecke algebra introduced in [C1]
and the quantum toroidal group,
as given in [GKV]. Indeed, one should keep in mind the geometrical
construction in [GRV] and [GKV] in terms of equivariant K-theory
of some flag manifolds. A similar K-theoretic construction of
Cherednik's algebra has motivated the present work.
At last, we would like to lay emphasis on the fact that,
contrary to [J], [D1] and [C2], the representations involved
in our duality are infinite dimensional.
Of course, in the classical case, i.e. $q=1$, a similar
duality holds between the toroidal Lie algebra and the toroidal
version of the symmetric group.

\vskip1cm

{\bf 1. Definition of the toroidal Hecke algebra.}

\vskip 3mm

For any non-negative integer $k$ set
$$[k]=\{0,1,2,...,k\}\qquad{\rm and}\qquad \set{k}=[k]\setminus\{0\}.$$

\vskip 3mm

{\bf 1.1. Definition.}
{\sl The toroidal Hecke algebra of type ${\got{gl}}\,(l)$, $\Htor$, is the
unital
associative algebra over $\Ac=\CC[\xb^{^{\pm 1}},\yb^{^{\pm 1}},\qb^{^{\pm
1}}]$
with generators
$$\ti^{^{\pm 1}},\quad\X_j^{^{\pm 1}},\quad \Y_j^{^{\pm 1}},\qquad
i\in\set{l-1}, \qquad j\in\set{l}$$
and the following relations :}
$$\ti\cdot\tii=\tii\cdot\ti=1,\qquad (\ti+1)\cdot (\ti-\qb^{^2})=0,$$
$$\ti\cdot\T_{i+1}\cdot\ti=\T_{i+1}\cdot\ti\cdot\T_{i+1},$$
$$\ti\cdot\T_j=\T_j\cdot\ti\quad {\sl if}\ |i-j|>1,$$
$$\X_0\cdot\Y_1=\xb\cdot\Y_1\cdot\X_0,
\qquad \X_i\cdot\X_j=\X_i\cdot\X_j,\qquad \Y_i\cdot\Y_j=\Y_j\cdot\Y_i,$$
$$\X_j\cdot\ti=\ti\cdot\X_j,\qquad \Y_j\cdot\ti=\ti\cdot\Y_j,
\qquad {\sl if}\ j\not= i,i+1$$
$$\ti\cdot\X_i\cdot\ti=\qb^{^2}\cdot\X_{i+1},\qquad
\tii\cdot\Y_i\cdot\tii=\qb^{^{-2}}\cdot\Y_{i+1},$$
$$\X_2\cdot\Y_1^{^{-1}}\cdot\X_2^{^{-1}}\cdot\Y_1=
\qb^{^{-2}}\cdot\yb\cdot\T_1^{^2},$$
{\sl where $\X_0=\X_1\cdot\X_2\cdots\X_l$.}

\vskip 3mm

{\bf 1.2.} \rems
{\bf (i)} When $\xb$ is taken to be 1 the toroidal Hecke algebra is nothing but
the double affine Hecke algebra introduced by Cherednik (see [C1]).

{\bf (ii)} Note that the map $\ti\mapsto\tii$,\quad $\X_i\mapsto\Y_i$,\quad
$\Y_i\mapsto\X_i$,\quad $\xb\mapsto\xb^{^{-1}}$,\quad $\yb\mapsto\yb^{^{-1}}$,
\quad $\qb\mapsto\qb^{^{-1}}$, extends to an automorphism, $\Xi$, of
$\Htor$ over $\CC$.

\vskip 3mm

{\bf 1.3.} Given $1\leq i \leq j <l$  set
$\T_{_{i,j}}=\T_i\cdot\T_{i+1}\cdots\T_j$ and
$\T_{_{j,i}}=\T_j\cdot\T_{j-1}\cdots\T_i$.
Then, put $\Q=\X_1\cdot\T_{_{1,l-1}}\in\Htor$. Clearly,
$\ti^{^{\pm 1}},\quad\Y_j^{^{\pm 1}},\quad \Q^{^{\pm 1}}$\quad
($i\in\set{l-1},\quad  j\in\set{l}$) is a system of generators of $\Htor$.
Besides, for any $i\in\set{l-1}$ a direct computation gives
$\Q\cdot\Y_i\cdot\Q^{^{-1}}=\yb^{^{-1}}\cdot\Y_{i+1}$, and
$\Q\cdot\Y_l\cdot\Q^{^{-1}}=\xb\cdot\yb^{^{l-1}}\cdot\Y_1$.
Indeed we have (see [C1])

\vskip.3cm

\prop {\sl The toroidal Hecke algebra $\Htor$ admits an
equivalent presentation in terms of generators
$$\ti^{^{\pm 1}},\quad\Y_j^{^{\pm 1}},\quad \Q^{^{\pm 1}},
\qquad  i\in\set{l-1},\qquad  j\in\set{l},$$
with relations :}
$$\ti\cdot\tii=\tii\cdot\ti=1,\qquad (\ti+1)\cdot (\ti-\qb^{^2})=0,$$
$$\ti\cdot\T_{i+1}\cdot\ti=\T_{i+1}\cdot\ti\cdot\T_{i+1},$$
$$\ti\cdot\T_j=\T_j\cdot\ti\quad {\sl if}\ |i-j|>1,$$
$$\Y_i\cdot\Y_j=\Y_j\cdot\Y_i,\qquad
\ti^{^{-1}}\cdot\Y_i\cdot\ti^{^{-1}}=\qb^{^{-2}}\cdot\Y_{i+1},$$
$$\Y_j\cdot\ti=\ti\cdot\Y_j,
\qquad {\sl if}\ j\not= i,i+1$$
$$\Q\cdot\T_{i-1}\cdot\Q^{^{-1}}=\ti\quad (1<i<l-1),\qquad
\Q^{^2}\cdot\T_{l-1}\cdot\Q^{^{-2}}=\T_1,$$
$$\Q\cdot\Y_i\cdot\Q^{^{-1}}=\yb^{^{-1}}\cdot\Y_{i+1}\quad (1\leq
i\leq l-1),\qquad  \Q\cdot\Y_l\cdot\Q^{^{-1}}=\xb\cdot\yb^{^{l-1}}\cdot\Y_1.$$

\vskip 3mm

{\bf 1.4.}  Let $\Hafa,\Hafb\i\Htor$ be the subalgebras generated respectively
by $\ti^{^{\pm 1}}$, $\Y_j^{^{\pm 1}}$  and
$\ti^{^{\pm 1}}$, $\X_j^{^{\pm 1}}$
($i\in\set{l-1}$, $j\in\set{l}$). Thus, $\Hafa$ and $\Hafb$ are
isomorphic to the affine Hecke algebra over $\Ac$
of type ${\got{gl}}\,(l)$, simply denoted $\Haf$.

\vskip 3mm

{\bf 1.5.}  Let $\Hf\i\Htor$ be the subalgebra generated by $\ti^{^{\pm 1}}$
($i\in\set{l-1}$). Thus, $\Hf$ is the Hecke algebra over $\Ac$ of finite type
${\got{gl}}\,(l)$.

\vskip 3mm

{\bf 1.6.}  Given complex numbers $q,x,y\in\CC^\times$,
let $\Mc_{_{x,y,q}}\i\Ac$ be
the maximal ideal generated by $\qb-q$, $\xb-x$ and $\yb-y$. Then, put
$\CC_{_{x,y,q}}=\Ac/\Mc_{_{x,y,q}}$ and
$$\htor=\Htor\otimes_{_\Ac}\CC_{_{x,y,q}},\qquad
\hf=\Hf\otimes_{_\Ac}\CC_{_{x,y,q}},$$
$$\hafa=\Hafa\otimes_{_\Ac}\CC_{_{x,y,q}},\qquad
\hafb=\Hafb\otimes_{_\Ac}\CC_{_{x,y,q}}.$$

\vskip 1cm

{\bf 2. Definition of the toroidal quantum group.}

\vskip 3mm

{\bf 2.1.} \df\ {\sl The toroidal quantum group of type ${\got{sl}}\,(n+1)$,
$\Utor$, is the unital associative algebra over
$\Bc=\CC[\cb^{^{\pm 1}},\db^{^{\pm 1}},\qb^{^{\pm 1}}]$
with generators}

$$\eb_{i,k},\quad\fb_{i,k},\quad\kb_{i,k},\quad\kb_{i,0}^{^{-1}},
\qquad i\in [n],\qquad k\in\ZZ.$$

{\sl The relations are expressed in term of the formal series}

$$\eb_i(z)=\sum_{k\in\interi}{\eb_{i,k}\cdot z^{^{-k}}},\quad
\fb_i(z)=\sum_{k\in\interi}{\fb_{i,k}\cdot z^{^{-k}}},\quad
\kb_i^{^\pm}(z)=\kb^{^{\pm 1}}_{i,0}+
\sum_{l\in\n}{\kb_{i,\pm l}\cdot z^{^{\mp l}}},$$

{\sl as follows}

$$\kb_{i,0}^{^{-1}}\cdot\kb_{i,0}=\kb_{i,0}\cdot\kb_{i,0}^{^{-1}}=1,\qquad
[\kb^{^\pm}_i(z),\kb^{^\pm}_j(w)]=0,\leqno(2.1.1)$$
$$\theta_{a_{ij}}(\cb^{^{-2}}\db^{^{m_{ij}}}\cdot z/w)\cdot\kb^{^+}_i(z)\cdot
\kb^{^-}_j(w)=\theta_{a_{ij}}(\cb^{^2}\db^{^{m_{ij}}}\cdot z/w)\cdot
\kb^{^-}_j(w)\cdot\kb^{^+}_i(z),\leqno(2.1.2)$$
$$\kb^{^\pm}_i(z)\cdot\eb_j(w) =\theta_{a_{ij}}
(\cb^{^{\pm 1}}\db^{^{m_{ij}}}\cdot z/w)\cdot\eb_j(w)\cdot\kb^{^\pm}_i(z),
\leqno(2.1.3)$$
$$\kb^{^\pm}_i(z)\cdot\fb_j(w) =\theta_{-a_{ij}}
(\cb^{^{\mp 1}}\db^{^{m_{ij}}}\cdot z/w)\cdot\fb_j(w)\cdot\kb^{^\pm}_i(z),
\leqno(2.1.4)$$
$$(\qb-\qb^{^{-1}})[\eb_i(z),\fb_j(w)]=\d_{ij}\,
\biggl(\d(\cb^{^{-2}}\cdot z/w)\cdot\kb^{^+}_i(\cb\cdot w)-
\d(\cb^{^2}\cdot z/w)\cdot\kb^{^-}_i(\cb\cdot z)\biggr),\leqno(2.1.5)$$
$$\eb_i(z)\cdot\eb_j(w)=\theta_{a_{ij}}
(\db^{^{m_{ij}}}z/w)\cdot\eb_j(w)\cdot\eb_i(z),\leqno(2.1.6)$$
$$\fb_i(z)\cdot\fb_j(w)=\theta_{-a_{ij}}
(\db^{^{m_{ij}}}z/w)\cdot\fb_j(w)\cdot\fb_i(z),\leqno(2.1.7)$$
$$\{(\eb_i(z_{_1})\cdot\eb_i(z_{_2})\cdot\eb_{j}(w)-
(\qb+\qb^{^{-1}})\cdot\eb_i(z_{_1})\cdot\eb_{j}(w)\cdot\eb_i(z_{_{2}})+
\eb_{j}(w)\cdot\eb_i(z_{_{1}})\cdot\eb_i(z_{_2})\}+\leqno(2.1.8)$$
$$+\{z_{_1}\leftrightarrow z_{_2}\}=0,
\qquad{\sl if}\quad |i-j|=1,$$
$$\{\fb_i(z_{_1})\cdot\fb_i(z_{_2})\cdot\fb_{j}(w)-
(\qb+\qb^{^{-1}})\cdot\fb_i(z_{_1})\cdot\fb_{j}(w)\cdot\fb_i(z_{_{2}})+
\fb_{j}(w)\cdot\fb_i(z_{_{1}})\cdot\fb_i(z_{_2})\}+\leqno(2.1.9)$$
$$+\{z_{_1}\leftrightarrow z_{_2}\}=0,
\qquad{\sl if}\quad |i-j|=1,$$

\vskip2mm

{\sl where $\d(z) = \sum_{_{n=-\infty}}^{^\infty}z^n$,
$\theta_m(z)={\qb^{^m}\cdot z -1\over z-\qb^{^m}}$ and
$a_{ij}$, $m_{ij}$, are the entries of the following $[n]\times[n]$-matrices }
\hfill\break

$\A=\pmatrix{
2&-1&0&&0&-1\cr
-1&2&-1&\cdots&0&0\cr
0&-1&2&&0&0\cr
&\vdots&&\ddots&&\vdots\cr
0&0&0&&2&-1\cr
-1&0&0&\cdots&-1&2\cr},
\qquad
\M=\pmatrix{
0&-1&0&&0&1\cr
1&0&-1&\cdots&0&0\cr
0&1&0&&0&0\cr
&\vdots&&\ddots&&\vdots\cr
0&0&0&&0&-1\cr
-1&0&0&\cdots&1&0\cr}.$
\hfill\break

\vskip 3mm

{\bf 2.2.} \rem
Note that, given $a\in\CC^\times$, the map
$\eb_{i,k}\mapsto a^{^k}\,\eb_{i+1,k},\quad
\fb_{i,k}\mapsto a^{^k}\,\fb_{i+1,k},\quad
\kb_{i,k}\mapsto a^{^k}\,\kb_{i+1,k}$\quad
$\eb_{n,k}\mapsto a^{^k}\,\eb_{0,k},\quad
\fb_{n,k}\mapsto a^{^k}\,\fb_{0,k},\quad
\kb_{n,k}\mapsto a^{^k}\,\kb_{0,k}$\quad
($i\in[n-1]$, $k\in\ZZ$),
extends to an automorphism, $\Psi_a$, of the $\Bc$-algebra $\Utor$.

\vskip 3mm

{\bf 2.3.}  Let $\Uafa,\Uafb\i\Utor$ be the subalgebras generated respectively
by $\eb_{i,k}$, $\fb_{i,k}$, $\kb_{i,k}$, $\kb_{i,0}^{^{-1}}$
($i\in\set{n}$, $k\in\ZZ$), and $\E_i=\eb_{i,0}$, $\F_i=\fb_{i,0}$,
$\K_i^{^{\pm 1}}=\kb_{i,0}^{^{\pm 1}}$ ($i\in[n]$).
The relations involving the Fourier coefficients of $\eb_i(z\,{\bf
d}^{^{-i}})$,
$\fb_i(z\,{\bf d}^{^{-i}})$ and $\kb^{^\pm}_i(z\,{\bf d}^{^{-i}})$\quad
$(i\in\set{n})$ are precisely the relations in
Drinfeld's ``new presentation" of the affine quantum group of type
${\widehat{\got{sl}}}\,(n+1)$ (see [D1]).
In another hand the relations involving $\eb_i$, $\fb_i$ and $\kb_i^{^\pm}$
are the quantum analogue of the Kac-Moody relations of type $A_n^{^{(1)}}$.
Thus $\Uafa$ and $\Uafb$ are both isomorphic to the affine quantum group of
type
${\widehat{\got{sl}}}\,(n+1)$, simply denoted by $\Uaf$ (see [B]).

\vskip 3mm

{\bf 2.4.}  Let $\Uf\i\Utor$ be the subalgebra generated by $\E_i$,
$\F_i$, $\K_i^{^{\pm 1}}$ ($i\in\set{n}$). Then $\Uf$ is
the quantum enveloping algebra of ${\got{sl}}\,(n+1)$.

\vskip 3mm

{\bf 2.5.} Given complex numbers $c,d,q\in\CC^\times$ let $\Nc_{_{c,d,q}}\i\Bc$
be the maximal ideal generated by $\cb-c$, $\db-d$, $\qb-q$.
Then, put $\CC_{_{c,d,q}}=\Bc/\Nc_{_{c,d,q}}$ and
$$\utor=\Utor\otimes_{_\Bc}\CC_{_{c,d,q}},\qquad
\uf=\Uf\otimes_{_\Bc}\CC_{_{c,d,q}},$$
$$\uafa=\Uafa\otimes_{_\Bc}\CC_{_{c,d,q}},\qquad
\uafb=\Uafb\otimes_{_\Bc}\CC_{_{c,d,q}}.$$
We say that a $\utor$-module $M$ is {\bf integrable} if
$$M=\bigoplus_{\l\in\ZZ^{^{[n]}}} M^{^\l},\qquad{\rm where}\quad
M^{^{(\l_0,\l_1,...,\l_n)}}=
\{m\in M\,|\,\kb_i\cdot m=q^{^{\l_i}}m,\quad\forall i\in[n]\},$$
and the $\eb_{i,k}$'s and $\fb_{i,k}$'s are locally nilpotent on $M$.
Moreover, a $\utor$-module is said to have
{\bf trivial central charge} if its restrictions
both to $\uafa$ and $\uafb$ have (in particular, $c=1$).

\vskip3mm

{\bf 2.6.} Fix $q\in\CC^\times$, and suppose that $l\leq n$.
Following [CP, 2.5], a $\uf$-module is said to be
{\bf of level} $l$ if its irreducible components are isomorphic to
irreducible components of $\V^{^{\ot l}}$. Then a $\uaf$-module or a
$\utor$-module is said to be of level $l$ if it is of level $l$ as a
$\uf$-module.

\vskip 1cm

{\bf 3. Definition of the duality functor.}

\vskip 3mm

Fix $c,d,x,q\in\CC^\times$ and set $y=1$.

\vskip 3mm

{\bf 3.1.} Let $\V$ be the fundamental
representation of $\uf$. It has a basis $\{v_{_i}\}_{_{i\in\set{n+1}}}$
on which the action of $\E_i,\F_i,\K_i$ ($i\in\set{n}$) is the following
$$\E_i( v_{_r})=\d_{r,i+1}\, v_{_{r-1}},
\qquad\F_i( v_{_r})=\d_{r,i}\, v_{_{r+1}},
\qquad\K_i( v_{_r})=q^{\d_{i,r}-\d_{i+1,r}}\, v_{_r}.$$
Then, $\V^{\ot l}$ is a left $\uf$-module for the induced action given by
the following coproduct
$$\Delta(\eb_i)=\eb_i\ot\,\kb_i+1\ot\,\eb_i,\qquad\Delta(\fb_i)=\fb_i\ot 1
+\kb^{^{-1}}_i\ot\,\fb_i,\qquad \Delta(\kb_i)=\kb_i\ot\,\kb_i.$$
This action commutes with the left $\hf$-action given by
$\ti=1^{^{\ot i-1}}\ot\,\T\,\ot 1^{^{\ot l-i-1}},$ where
$\T\in {\rm End}\,\V^{^{\ot 2}}$
verifies
$$\T(v_{_r}\ot v_{_s})=\cases{q^{^2}v_{_r}\ot v_{_s}&if $r=s,$\cr
q v_{_s}\ot v_{_r}&if $r<s,$\cr
q v_{_s}\ot v_{_r}+(q^{^2}-1)v_{_r}\ot v_{_s}&if $r>s.$\cr}$$

\vskip 3mm

{\bf 3.2.} For any $i\in\set{n+1}$ define $\tb_i\pr$ to be the
automorphism of the algebra $\uaf$ given on the Kac-Moody
generators by the formula
$$\tb_i\pr(\eb_i)=-\fb_i\cdot\kb_i,\qquad
\tb_i\pr(\eb_j)=\sum_{s=0}^{-a_{ij}}(-1)^{^{s-a_{ij}}}q^{^{-s}}
\eb_i^{^{(-a_{ij}-s)}}\cdot\eb_j\cdot\eb_i^{^{(s)}}\quad (i\not=j),$$
$$\tb_i\pr(\fb_i)=-\kb_i^{^{-1}}\cdot\eb_i,\qquad
\tb_i\pr(\fb_j)=\sum_{s=0}^{-a_{ij}}(-1)^{^{s-a_{ij}}}q^{^{s}}
\fb_i^{^{(s)}}\cdot\fb_j\cdot\fb_i^{^{(-a_{ij}-s)}}\quad (i\not=j),$$
$$\tb_i\pr(\kb_j)=\kb_{s_i(j)},$$
where
$\eb_i^{^{(j)}},
\fb_i^{^{(j)}}$ are the usual quantum divided powers (see [L1, 3.1.1])
and $s_i\in{\got{S}}_{n+1}$ is the transposition $(i\ i+1)$.
Let $M\pr$ be an integrable $\uaf$-module.
Set $\tb\ppr_i\in{\rm Aut}_{_{\CC}}(M\pr)$
($i\in\set{n+1}$) to be the braid operator defined by
$$\tb_i\ppr(m\pr)=\sum_{_{r-s+t=-k}}
(-1)^{^{s+k}}q^{^{s-rt}}\eb_i^{^{(r)}}\cdot
\fb_i^{^{(s)}}\cdot\eb^{^{(t)}}_i\cdot m',$$
where $m\pr\in M\pr$ and $k\in\NN$
are such that $\kb_i\cdot m\pr=q^{^{k}}m\pr$
(see [L1, 5.2.1 and  5.2.3]).
We have (see [L1, chapters 5 and 37])
$$\forall m'\in M',\quad\forall u\in\uaf,\qquad
\tb_i\ppr(u\cdot m')=\tb_i\pr(u)\cdot\tb_i\ppr(m').\leqno(3.2.1)$$

\vskip2mm

Similarly, denote by $\tau$ the automorphism of the affine Dynkin diagram
$A_n^{^{(1)}}$ (with vertices indexed by $1\vir n,n+1$) given
by
$\tau(i)=i+1\ (i\in\set{n}),\ \tau(n+1)=1.$
Let $\tau\pr$ be the automorphism of the algebra $\uaf$ given,
in terms of its Kac-Moody generators $\eb_i,\fb_i,\kb_i,\
(i\in\set{n+1})$, by the following rule
$$\tau\pr(\eb_i)=\eb_{\tau(i)},\quad
\tau\pr(\fb_{i})=\fb_{\tau(i)},\quad \tau\pr(\kb_i)=\kb_{\tau(i)}.$$
Put
$\tb\pr_{\o_1}=\tau\pr\circ\tb\pr_n\circ\tb\pr_{n-1}\circ
\cdots\circ\tb\pr_1\in{\rm Aut}
(\uaf).$
Take a right $\haf$-module $M$. In particular, $M$ is a right $\hf$-module
and we can consider the dual left $\uf$-module $\modulo.$
This module is endowed with a structure
of left $\uaf$-module such that for any $m\in M$
and $\vu\in \V^{^{\ot l}}$ (see [CP])
$$\eb_{n+1}(m\ot\vu)=\sum_{j=1}^lm\cdot\Y^{^{-1}}_j\ot\,\fb_{_{\t,j}}(\vu),
\qquad\fb_{n+1}(m\ot \vu)=\sum_{j=1}^lm\cdot\Y_j\ot\,\eb_{_{\t,j}}(\vu),$$
$$\kb_{n+1}(m\ot\vu)=m\ot\,(\kt^{^{-1}})^{^{\ot l}}(\vu),$$
where $\et,\ft,\kt\in{\rm End}_{_\CC}(\V)$ are defined by means of
$$\et\cdot v_r=\d_{r,n+1}\,v_1,\qquad \ft\cdot v_r=\d_{r,1}\,v_{n+1},
\qquad\kt\cdot v_r=q^{\d_{1,r}-\d_{n+1,r}}\,v_r,$$
and  $\fb_{_{\t,j}}=1^{^{\ot j-1}}\ot\,\ft\ot\,(\kt^{^{-1}})^{^{\ot l-j}}$,
$\eb_{_{\t,j}}=\kt^{^{\ot j-1}}\ot\,\et\ot 1^{^{\ot l-j}}$.
Until the end of section 3.2 take $M\pr=\modulo.$
Define
$\tau\ppr\in{\rm Aut} (M\pr)$ such that
$$\tau\ppr(m\ot v_\jb)=m\cdot \Y_1^{^{\d_{n+1,j_1}}}\cdot
\Y_2^{^{\d_{n+1,j_2}}}\cdots\Y_l^{^{\d_{n+1,j_l}}}\ot
v_{j_1+1}\ot v_{j_2+1}\ot\cdots\ot v_{j_l+1}$$
where $v_{n+2}$ stands for $v_1$ and $\jb=(j_1\vir j_l)$
is an $l$-tuple of integers in $\set{n+1}.$
A direct computation gives
$$\forall m\pr\in M\pr,\quad\forall u\in\uaf,\qquad
\tau\ppr(u\cdot m\pr)=\tau\pr(u)\cdot\tau\ppr(m\pr).\leqno(3.2.2)$$
Put $\tb_{\o_1}\ppr=\tau\ppr\circ\tb_n\ppr\circ\tb_{n-1}\ppr\circ
\cdots\circ\tb_1\ppr \in{\rm Aut}(M\pr).$
As a consequence of (3.2.1) and (3.2.2),
$$\forall m\pr\in M\pr,\quad\forall u\in\uaf,\qquad
\tb\ppr_{\o_1}(u\cdot m\pr)=\tb_{\o_1}\pr(u)\cdot\tb_{\o_1}\ppr
(m\pr).\leqno(3.2.3)$$

{\bf Example.}
If $l=1$ we find
$$\tb\ppr_i(m\ot v_j)=(-1)^{^{\d_{i+1,j}}}q^{^{\d_{ij}}}m\ot v_{s_i(j)}\qquad
i\in\set{n},\quad j\in\set{n+1},\leqno(3.2.4)
$$
$$\tb\ppr_{\o_1}(m\ot v_j)= -m\cdot (-q^{^{n}}\Y_1)^{^{\d_{1,j}}}\ot v_j.
\leqno(3.2.5)$$

For any $l$-tuple $\jb\in(\set{n+1})^l$ set
$\vu_{\jb}=v_{_{j_1}}\ot\cdots\ot v_{_{j_l}}.$

\vskip 3mm

\lem {\sl Let $\jb$ be non decreasing,  and set
$ \Y_{1,s}=\Y_1\cdot\Y_2\cdots\Y_s$, with
$\jb^{^{-1}}(1)=]0,s].$ Then }
$$\tb^{^{\pp\pp }}_{\o_1}(m\ot\vb_\jb)=
(-1)^{^{l+s}}q^{^{ns}}m\cdot\Y_{_{1,s}}\ot\vb_\jb.
$$
\vskip 3mm
\pf
{}From [L1, 5.3.4], for all integrable $\uaf$-modules $N_1,N_2$ and
for all $n_1\in N_1,n_2\in N_2$
$$ \eb_i(n_1)\ot\fb_i(n_2)=0\ \Longrightarrow\  \tb_i\ppr(n_1\ot n_2)
=\tb_i\ppr(n_1)\ot\tb_i\ppr(n_2)\leqno(3.2.6)$$
(note that Lusztig uses the  opposite coproduct in [L1]).
Put $\jbp=(j_1\vir j_{l-1})$ and for any permutation $\si\in
{\got{ S}}_{n+1}$, set $\si(\jb)=(\si(j_1)\vir \si(j_l)).$
We first compute
$\tb_n\ppr\circ\tb_{n-1}\ppr\circ\cdots\circ\tb_1\ppr$.
Note that
$\tb_n\ppr\circ\tb_{n-1}\ppr\circ\cdots\circ\tb_1\ppr (m\ot\vu_\jb)
=m\ot
\tb_n\ppr\circ\tb_{n-1}\ppr\circ\cdots\circ\tb_1\ppr (\vu_\jb).$
We have  $\eb_1(\vu_\jbp)\not=0$ if and only if
$\jbp^{^{-1}}(2)\not=\emptyset$ and  $\fb_1(v_{j_l})\not=0$
if and only if $j_l=1$.
Since $\jb$ is non decreasing, using (3.2.6) and (3.2.4) we find
$$\tb_1\ppr(\vu_\jb)=\tb_1\ppr(\vu_\jbp)\ot\tb_1\ppr(v_{j_l})=
\cdots=
\bigotimes_{k=1}^l\tb_1\ppr(v_{j_k})=(-1)^{^{a_2}}q^{^{a_1}}\vu_{s_1(\jb)},$$
where $a_i$ is the cardinality of $\jb^{^{-1}}(i)$.
Suppose  that
$$\tb\ppr_k\circ\cdots\circ\tb\ppr_1(\vu_\jb)=
(-1)^{^{a_2+\cdots+a_{k+1}}}q^{^{ka_1}}\vu_{s_k\cdots
s_1(\jb)}.$$
Now
$\eb_{k+1}(\vu_{s_k\cdots s_1(\jbp)})\ot
\fb_{k+1}(v_{s_k\cdots s_1(j_l)})
\not=0$
if and only if
$(s_k\cdots s_1(\jbp))^{^{-1}}(k+2)\not=\emptyset$ and $j_l=1,$
that is if and only if $\jbp^{^{-1}}(k+2)\not=\emptyset$ and $j_l=1$.
Since $\jb$ is non decreasing, using (3.2.6) and (3.2.4)
we get the same formula as above for
$\tb\ppr_{k+1}\circ\cdots\circ\tb\ppr_1(\vu_\jb).$
Then, finally
$$\tb\ppr_n\circ\cdots\circ\tb\ppr_1(m\ot\vu_\jb)=
(-1)^{^{l-s}}q^{^{ns}}m\ot v_{_{j_1-1}}\ot\cdots\ot v_{_{j_l-1}},  $$
where $v_{_0}$ stands for $v_{_{n+1}}$.
Applying $\tau\ppr$ we find the
result.
\cvd

\vskip 3mm

{\bf 3.3.}
Take now a right $\hafa$-module $M$ and consider its
dual left $\uafa$-module $\modulo$. Since $\uafa\simeq\uaf$ (see 2.3)
we can use all the constructions of the section 3.2 for the
$\uafa$-module $\modulo.$ In particular, the formulas for $\eb_{n+1},$
$\fb_{n+1}$ and $\kb_{n+1}$ write as follows
$$\eb_{n+1}(m\ot\vu)=d^{^{-1}}\sum_{j=1}^lm\cdot\Y^{^{-1}}_j\ot\,
\fb_{_{\t,j}}(\vu),
\qquad\fb_{n+1}(m\ot \vu)=d\sum_{j=1}^lm\cdot\Y_j\ot\,\eb_{_{\t,j}}(\vu),$$
$$\kb_{n+1}(m\ot\vu)=m\ot\,(\kt^{^{-1}})^{^{\ot l}}(\vu).$$

The $\uafa$-module structure on $\modulo$ is given in terms
of Drinfeld generators by (see [GRV]) :

\vskip 3mm

\th {\sl If $\jb$ is non decreasing and $i\in\set{n}$ we have
$$\matrix{
{\dd \eb_i(z)(m\ot \vu_\jb)=q^{^{1-t+s}}m\cdot
(1+\sum_{k=s+1}^{t-1}\T_{_{k,s+1}})\cdot
\d(q^{^{n+1-i}}d^{^i}z\Y_{s+1})\ot\vu_{\jb_{s+1}^-},}\hfill\cr\cr
{\dd \fb_i(z)(m\ot \vu_\jb)=q^{^{1-s+r}}m\cdot
(1+\sum_{k=r+1}^{s-1}\T_{_{k,s-1}})\cdot
\d(q^{^{n+1-i}}d^{^i}z\Y_{s})\ot\vu_{\jb_{s}^+},}\hfill\cr\cr
{\dd \kb^{^\pm}_i(z)(m\ot\vu_\jb)=m\cdot\prod_{j_k=i}{
\theta^{^\pm}_{1}(q^{^{n+2-i}}d^{^i}z\Y_k)\cdot
\prod_{j_k=i+1}{\theta^{^\pm}_{-1}(q^{^{n-i}}}d^{^i}z\Y_k)}\ot\vu_\jb,}\hfill\cr
}$$
where $]r,s]=\jb^{-1}(i),\quad]s,t]=\jb^{-1}(i+1)$,
$$\cases{\vu_{\jb_{s+1}^-}=0&if $\quad \jb^{-1}(i+1)=\emptyset,$\cr
\jb_{s+1}^-=(j_1\vir j_s,i,j_{s+2}\vir j_l)& else,\cr}$$
$$\cases{\vu_{\jb_{s}^+}=0&if $\quad \jb^{-1}(i)=\emptyset,$\cr
\jb_{s}^+=(j_1\vir j_{s-1},i+1,j_{s+1}\vir j_l)& else,\cr}$$
and
$\t^{^\pm}_m$ stands for the Taylor expansion of $\t_m$
respectively at $\infty$ and $0$.}

\vskip 3mm

\pf  Let us  prove  the formula involving $\eb_{i,h}$.
We know that $\eb_{1,h}=(-d)^{^{-h}}\tb_{\o_1}^{^{\pp -h}}(\eb_1)$
(see [B, 4.6] and 2.3). So (3.2.3) gives
$$\eb_{1,h}\circ\,\tb_{\o_1}^{^{\pp\pp -h}}(m\ot\vu_\jb)
=(-d)^{^{-h}}\tb_{\o_1}^{^{\pp\pp -h}}\circ\,
\eb_1(m\ot\vu_\jb).\leqno(3.3.1)$$
First of all note that
formulas in section 3.1 give  (with $\jb^{^{-1}}(1)=]0,s]$ and
$\jb^{^{-1}}(2)=]s,t]$)
$$\eb_1(m\ot \vu_\jb)=q^{^{1-t+s}}m\cdot
(1+\sum_{k=s+1}^{t-1}\T_{_{k,s+1}})\ot\vu_{\jb_{s+1}^-}.$$

We have
$$\tb_{\o_1}^{^{\pp\pp -h}}\circ\,\eb_1\,(m\ot \vu_\jb)=
q^{^{1-t+s}}\tb_{\o_1}^{^{\pp\pp -h}}\circ
(1+\sum_{k=s+1}^{t-1}\T_{_{k,s+1}})\,(m\ot\vu_{\jb^-_{s+1}})=$$
$$=q^{^{1-t+s}}(1+\sum_{k=s+1}^{t-1}\T_{_{k,s+1}})\circ
\tb_{\o_1}^{^{\pp\pp -h}} (m\ot\vu_{\jb^-_{s+1}}),$$
since  $\T_i$ and $\tb\ppr_{\o_1}$ commute.
Thus, using lemma 3.2,  the right hand side of (3.3.1) is
$$(-1)^{^{h(s+1+l)}}(-d)^{^{-h}}q^{^{1-t+s-h(s+1)n}}
(1+\sum_{k=s+1}^{t-1}\T_{_{k,s+1}})\,
(m\cdot\Y_{_{1,s+1}}^{^{-h}}\ot\vu_{\jb_{s+1}^-}).$$
A similar computation for the left hand side of (3.3.1) gives
$$(-1)^{^{h(s+l)}}q^{^{-hsn}}\eb_{1,h}\,
(m\cdot \Y_{_{1,s}}^{^{-h}}\ot\vu_\jb)=
(-1)^{^{h(s+l)}}q^{^{-hsn}}\Y_{_{1,s}}^{^{-h}}
\cdot\eb_{1,h}\,(m\ot\vu_\jb),$$
since $\Y_i$ commutes to $\eb_{1,h}$.
Now $\Y_i$ commutes also with $\T_{_{k,s+1}}$ if $i\in\set{s}$.
Thus we finally obtain
$$\eb_{1,h}\,(m\ot\vu_\jb)=
d^{^{-h}}q^{^{1-t+s-hn}}(1+\sum_{k=s+1}^{t-1}\T_{_{k,s+1}})
\cdot\Y_{s+1}^{^{-h}}\,(m\ot\vu_{\jb^-_{s+1}}).$$
For the other cases we  proceed in a similar way. Namely,
since  the $i$-th fundamental weight $\o_i$
 verifies $\o_i=\tau^{^{i}}\cdot(s_n\cdot s_{n-1}\cdots
s_1)^{^{i}}$, define $\tau_{\o_i}\ppr=\tau^{^{\pp\pp i}}\circ
(\tb_n\ppr\circ\tb_{n-1}\ppr\circ\cdots\circ\tb_1\ppr)^{^{i}}.$
Then $\eb_{i,h}=(-d)^{^{-ih}}\tb_{\o_i}^{^{\pp\pp -h}} (\eb_i).$
\cvd

\vskip 3mm

\rem The relations (2.1.3) give in particular
$$\eb_1\cdot\kb_{2,1}-q\,\kb_{2,1}\cdot\eb_1=
(q-q^{^{-1}})c^{^{-1}}d^{^{-1}}\eb_{1,1}\cdot\kb_2.$$
Moreover, for any non-decreasing $l$-tuple $\jb$ one gets
$$\kb_i(m\otimes\vb_\jb)=q^{^{\sharp\jb^{^{-1}}(i)-\sharp\jb^{-1}(i+1)}}
m\otimes\vb_\jb,$$
$$\kb_{i,1}(m\otimes\vb_\jb)=q^{^{i-n+\sharp\jb^{-1}(i)-\sharp\jb^{-1}(i+1)}}
(1-q^{^{-2}})d^{^{-i}}m\cdot(q^{^{-1}}\sum_{_{\jb_k=i}}\Y_k^{^{-1}}-
q\sum_{_{\jb_k=i+1}}\Y_k^{^{-1}})\otimes\vb_\jb,$$
Thus, the evaluation of the above formula on a non-zero vector of type
$m\otimes v_1\otimes v_2\otimes\cdots\otimes v_l$ gives $c=1$.
In other words, on the dual module $\modulo$ the central element $c$ is
trivial, even if $M$ is not finite dimensional.

\vskip 3mm

{\bf 3.4.} In the two next sections set $c=1$ and
take a right $\htor$-module $M$. In particular $M$ is a right
$\hafa$-module and thus $\modulo$ is endowed with a $\uafa$-module structure.
Let $\psi\,:\,\modulo\rightarrow\modulo$ be the linear map
defined for any $l$-tuple $\jb=(j_1,j_2,...,j_l)$ and any $m\in M$ by
$$\psi(m\otimes\vb_\jb)=m\cdot\X_1^{^{-\d_{n+1,j_1}}}
\cdot\X_2^{^{-\d_{n+1,j_2}}}\cdots\X_l^{^{-\d_{n+1,j_l}}}
\otimes v_{_{1+j_1}} \otimes v_{_{1+j_2}}\otimes\cdots\otimes v_{_{1+j_l}}$$
where $v_{_{n+2}}$ stands for $v_{_1}$.

\vskip3mm

\prop {\sl Given $i\in[n]\setminus\{0,1\}$, we have the following identities
in ${\rm End}(\modulo)$}
$$\matrix{
\psi^{^{-1}}\circ\eb_i(z)\circ\psi = \eb_{i-1}(q^{^{-1}}dz),\qquad\hfill
&\psi^{^{-2}}\circ\eb_1(xz)\circ\psi^{^2} =
\eb_{n}(q^{^{n-1}}d^{^{1-n}}z),\hfill\cr
\psi^{^{-1}}\circ\fb_i(z)\circ\psi = \fb_{i-1}(q^{^{-1}}dz),\qquad\hfill
&\psi^{^{-2}}\circ\fb_1(xz)\circ\psi^{^2} =
\fb_{n}(q^{^{n-1}}d^{^{1-n}}z),\hfill\cr
\psi^{^{-1}}\circ\kb_i(z)\circ\psi = \kb_{i-1}(q^{^{-1}}dz),\qquad\hfill
&\psi^{^{-2}}\circ\kb_1(xz)\circ\psi^{^2} =
\kb_{n}(q^{^{n-1}}d^{^{1-n}}z).\hfill\cr
}$$

\vskip3mm

In order to prove this proposition we need the following

\vskip 3mm

{\bf  Lemma.} {\sl For any $1\leq i\leq j<l$ put
$\Q_{_{i,j}}=\X_{i}\cdot\T_{_{i,j}}\in\Htor$.
Then, if $i\leq r\leq j$ and $i<t<j$
$$\Q_{_{i,j}}\cdot \Y_r\cdot \Q_{_{i,j}}^{^{-1}}=\yb^{^{-1}}\cdot\Y_{r+1},
\qquad \Q_{_{i,j}}\cdot\T_{t-1}\cdot\Q_{_{i,j}}^{^{-1}}=\T_t.$$}

\vskip 1mm

\pf We prove  the first equality, the second being similar.
First, suppose that $i=1$ and use a decreasing induction on $j$.
If $j=l-1$ we get the relation
$\Q\cdot\Y_r\cdot\Q^{^{-1}}=\yb^{^{-1}}\cdot\Y_{r+1}.$
Take $r$ such that  $1\leq r\leq j-1$, then $r\leq j$ and by induction we get
$$\X_1\cdot\T_{_{1,j-1}}\cdot\T_j\cdot\Y_r\cdot\T_j^{^{-1}}\cdot
\T_{_{1,j-1}}^{^{-1}}\cdot\X_1^{^{-1}}=\yb^{^{-1}}\cdot\Y_{r+1}.$$
Since $r\neq j,j+1$, we have $\T_j\cdot\Y_r=\Y_r\cdot\T_j$ and we are done.

Fix now $j$ and make induction on $i$, the case $i=1$ being proved before.
Consider $i+1\leq r\leq j$. Then $i\leq r$ and by induction we have
$$\X_i\cdot\ti\cdot\T_{_{i+1,j}}\cdot\Y_r\cdot\T_{_{i+1,j}}^{^{-1}}\cdot
\tii\cdot\X_i^{^{-1}}=\yb^{^{-1}}\cdot\Y_{r+1}.$$
Now $\X_i\cdot\ti=q^{^{2}}\tii\cdot\X_{i+1},$ then we find
$$\tii\cdot (\X_{i+1}\cdot\T_{_{i+1,j}}\cdot\Y_r\cdot\T_{_{i+1,j}}^{^{-1}}
\cdot\X_{i+1})\cdot\ti=\yb^{^{-1}}\cdot\Y_{r+1},$$
i.e.,
$$\Q_{_{i+1,j}}\cdot\Y_r\cdot\Q_{_{i+1,j}}^{^{-1}}=\yb^{^{-1}}\cdot\ti\cdot
\Y_{r+1}\cdot\tii.$$
Since $1+r\neq i,i+1$ we have $\ti\cdot\Y_{r+1}=\Y_{r+1}\cdot\ti$ and we are
through.\cvd

\vskip 3mm

{\it Proof of the proposition.} We prove  the relation involving $\eb_i$.
Suppose first that $i\neq 0,1$. Take a non decreasing $l$-tuple $\jb$. Put
$$\jb^{^{-1}}(i-1)=]r,s],\qquad
\jb^{^{-1}}(i)=]s,t],\qquad\jb^{^{-1}}(n+1)=]p,l],$$
$$\jbp=(j_1+1,j_2+1\vir j_p+1,1\vir 1),\qquad
\jbpp= (1\vir 1,j_1+1,j_2+1\vir j_p+1).$$
Thus $\jbpp^{^{-1}}(i)=]l-p+r,l-p+s]$ and $\jbpp^{^{-1}}(i+1)=]l-p+s,l-p+t]$.
Set $\R_p=q^{^{p(p-l)}}\T_{_{p,1}}\cdot\T_{_{p+1,2}}\cdots\T_{_{l-1,l-p}}$.
Then
$$\eb_i(z)\circ\psi(m\ot\vu_\jb)=\eb_i(z)
(m\cdot\X_{p+1}^{^{-1}}\cdot\X_{p+2}^{^{-1}}\cdots\X_l^{^{-1}}\ot\vu_\jbp)=
\eb_i(z)(m\cdot\X_{p+1}^{^{-1}}\cdot\X_{p+2}^{^{-1}}\cdots\X_l^{^{-1}}
\cdot\R_p\ot\vu_\jbpp)=$$
$$=q^{^{1-t+s}}m\cdot\X_{p+1}^{^{-1}}\cdots\X_l^{^{-1}}\cdot\R_p\cdot
(1+\sum_{k=l-p+s+1}^{l-p+t-1}\T_{_{k,l-p+s+1}})
\cdot\d(q^{^{n+1-i}}d^{^{i}}z\Y_{l-p+s+1})\ot\vu_\jbppp,$$
where $\jbppp=(\jbpp)^-_{l-p+s+1}.$ On the other side
$$\psi\circ\eb_{i-1}(q^{^{-1}}dz)(m\ot\vu_\jb)=\psi\biggl(q^{^{1-t+s}}m\cdot
(1+\sum_{k=s+1}^{t-1}\T_{_{k,s+1}})\cdot
\d(q^{^{n+1-i}}d^{^{i}}z\Y_{s+1}\ot\vu_{\jb_{s+1}^-})\biggr)=$$
$$=q^{^{1-t+s}}m\cdot(1+\sum_{k=s+1}^{t-1}\T_{_{k,s+1}})
\cdot\d(q^{^{n+1-i}}d^{^{i}}z\Y_{s+1})\cdot
\X_{p+1}^{^{-1}}\cdots\X_l^{^{-1}}\cdot\R_p\ot\vu_\jbppp.$$
Now, $\R^{^{-1}}_p\cdot\X_l\cdot\X_{l-1}\cdots\X_{p+1}=
q^{^{p(p-l)}}\P_p$,
where $\P_p=\Q_{_{l-p,l-1}}\cdots\Q_{_{2,p+1}}\cdot\Q_{_{1,p}}$.
Thus the relation follows from
$$\P_p\cdot(1+\sum_{k=s+1}^{t-1}\T_{_{k,s+1}})\cdot\d(\Y_{s+1})\cdot
\P_p^{^{-1}}=
(1+\sum_{k=l-p+s+1}^{l-p+t-1}\T_{_{k,l-p+s+1}})\cdot\d(\Y_{l-p+s+1}),
$$
which is a consequence of preceding lemma (note that $1\leq s+1\leq p$
and  $0<k<p-1$ for any $k\in ]s,t[$).

Suppose now that $i=1$. Set
$$\jb^{^{-1}}(n)=]r,s],\qquad\jb^{^{-1}}(n+1)=]s,l],$$
$$\jbp=(1\vir 1,2\vir 2, j_1+2\vir j_r+2),\qquad\jbp^{^{-1}}(1)=]0,s-r],\qquad
\jbp^{^{-1}}(2)=]s-r,l-r],$$
$$\jbpp=(1\vir 1,2\vir 2,j_1+2\vir j_r+2),\quad \jbpp^{^{-1}}(1)=
]0,s-r+1],\quad \jbpp^{^{-1}}(2)=]s-r+1,l-r],$$
$$\jbppp=(j_1\vir j_r,n\vir n,n+1\vir n+1),\qquad \jbppp^{^{-1}}(n)=]r,s+1],
\qquad \jbppp^{^{-1}}(n+1)=]s+1,l],$$
$$\tbp=(j_1+2\vir j_r+2,1\vir 1,2\vir 2),\qquad \tbp^{^{-1}}(1)=]r,s+1],
\qquad \tbp^{^{-1}}(2)=]s+1,l].$$
Then
$$\eb_1(xz)\circ\psi^{^{2}}(m\ot\vu_\jb)=
\eb_1(xz)(m\cdot\X_{r+1}^{^{-1}}\cdots\X_l^{^{-1}}\cdot\R_r\ot\vu_\jbp)=$$
$$=q^{^{1-l+s}}m\cdot\X_{r+1}^{^{-1}}\cdots\X_l^{^{-1}}\cdot\R_r\cdot(1+
\sum_{k=s-r+1}^{l-r-1}\T_{_{k,s-r+1}})\cdot\d(q^{^{n}}dxz
\Y_{s-r+1})\ot\vu_\jbpp.$$
In another hand,
$$\psi^{^{2}}\circ\eb_n(q^{^{n-1}}d^{^{1-n}}z)(m\ot\vu_\jb)
=\psi^{^{2}}
\biggl(q^{^{1-l+s}}m\cdot(1+\sum_{k=s+1}^{l-1}\T_{_{k,s+1}})
\cdot\d(q^{^{n}}dz\Y_{s+1})\ot\vu_\jbppp\biggr)=$$
$$=q^{^{1-l+s}}m\cdot(1+\sum_{k=s+1}^{l-1}\T_{_{k,s+1}})\cdot\d(q^{^{n}}dz\Y_{s+1})
\cdot\X_{r+1}^{^{-1}}\cdots\X_l^{^{-1}}\ot\vu_\tbp=$$
$$=q^{^{1-l+s}}m\cdot(1+\sum_{k=s+1}^{l-1}\T_{_{k,s+1}})\cdot\d(q^{^{n}}dz\Y_{s+1})
\cdot\X_{r+1}^{^{-1}}\cdots\X_l^{^{-1}}\cdot\R_r\ot\vu_\jbpp.$$
Thus it suffices to prove that
$$\P_r\cdot(1+\sum_{k=s+1}^{l-1}\T_{_{k,s+1}})\cdot\d(\Y_{s+1})\cdot\P_r^{^{-1}}=
(1+\sum_{k=s-r+1}^{l-r-1}\T_{_{k,s-r+1}})\cdot\d(x\Y_{s-r+1}).
\leqno{(3.4.1)}$$
Formula (3.4.1) follows from

\vskip 3mm

{\bf Lemma.} {\sl For any $r<l$ set
$\P_r=\Q_{_{l-r,l-1}}\cdots\Q_{_{2,r+1}}\cdot\Q_{_{1,r}}\in\Htor$. Then, if\
$r<s+1$ and $r<k<l$,}
$$\P_r\cdot\Y_{s+1}\cdot\P_r^{^{-1}}=\xb\cdot\yb^{^{r}}\cdot\Y_{s-r+1},\qquad
\P_r\cdot\T_k\cdot\P_r^{^{-1}}=\T_{k-r}.$$

\vskip 1mm

\pf By definition of $\P_r$ one gets
$$\X_0\cdot\P_r^{^{-1}}=q^{^{2r(l-r)}}\X_1\cdot\X_2\cdots\X_r\cdot
\T_{_{r,1}}\cdot\T_{_{r+1,2}}\cdots\T_{_{l-1,l-r}}.$$
Then a direct computation gives
$\P_r=q^{^{2r(r-l)}}\Q_{_{1,l-r}}^{^{-1}}\cdot\Q_{_{2,l-r+1}}^{^{-1}}\cdots
\Q_{_{r,l-1}}^{^{-1}}\cdot\X_0,$
and the result follows from the preceding lemma.
\cvd

\vskip 3mm

{\bf 3.5.} Take $\eb_0(z),\fb_0(z),\kb_0(z)\in
{\rm End}(\modulo)[[z^{^{\pm 1}}]]$ such that
$$\eb_0(z)=\psi^{^{-1}}\circ\eb_1(qd^{^{-1}}z)\circ\psi,\quad
\fb_0(z)=\psi^{^{-1}}\circ\fb_1(qd^{^{-1}}z)\circ\psi,\quad
\kb_0(z)=\psi^{^{-1}}\circ\kb_1(qd^{^{-1}}z)\circ\psi.$$

The operators are thus defined in such a way that if $i\in[n]$
and $x=d^{^{-n-1}}q^{^{n+1}}$,
$$\psi^{^{-1}}\circ\eb_i(z)\circ\psi = \eb_{i-1}(q^{^{-1}}dz),\qquad
\psi^{^{-1}}\circ\fb_i(z)\circ\psi = \fb_{i-1}(q^{^{-1}}dz),
\leqno(3.5.1)$$
$$\psi^{^{-1}}\circ\kb_i^{^\pm}(z)\circ\psi = \kb_{i-1}^{^\pm}(q^{^{-1}}dz),$$

where $\eb_{-1}(z)$, $\fb_{-1}(z)$ and $\kb^{^\pm}_{-1}(z)$ stand for
$\eb_n(z)$, $\fb_n(z)$ and $\kb^{^\pm}_n(z)$.

\vskip 3mm

{\bf Remark.}
Let  $\et,\ft,\kt\in{\rm End}_{_\CC}(\V)$ and
 $\fb_{_{\t,j}},
\eb_{_{\t,j}}$ be defined as in 3.2.
A direct computation gives, for any $\vu\in\V^{^{\ot l}}$ and $m\in M$,
$$\eb_0(m\ot\vu)=\sum_{j=1}^lm\cdot\X_j\ot\,\fb_{_{\t,j}}(\vu),
\qquad\fb_0(m\ot \vu)=\sum_{j=1}^lm\cdot\X_j^{^{-1}}\ot\,\eb_{_{\t,j}}(\vu),$$
$$\kb_0(m\ot\vu)=m\ot\,(\kt^{^{-1}})^{^{\ot l}}(\vu).$$
In particular $\kb_0\cdot\kb_1\cdots\kb_n(m\ot\vu)=m\ot\vu.$

\vskip 3mm

The main result of this section is the following theorem :

\vskip3mm

\th {\sl Suppose that $x=d^{^{-n-1}}q^{^{n+1}}$ and $c=y=1$.
Then for any right $\htor$-module $M$, the preceding formulas give
a left integrable $\utor$-module structure on $\modulo$.
Moreover $\modulo$ has trivial central charge.}

\vskip 3mm

{\sl Proof.} By construction, the formulas for
$\eb_i(z),\fb_i(z),\kb_i^{^\pm}(z)\in{\rm End}(\modulo)[[z^{^{\pm 1}}]],$
$i\in\set{n}$, written in theorem 3.3 give a representation of $\uafa.$
In order to verify all the relations 2.1, just notice that it's
sufficient to verify the relations involving the generators
$\Psi^{^k}_{_{qd^{^{-1}}}}(\eb_i(z)),$ $\Psi^{^k}_{_{qd^{^{-1}}}}(\fb_i(z))$
and $\Psi^{^k}_{_{qd^{^{-1}}}}(\kb_i^{^{\pm}}(z))$
of the subalgebras
$\Psi^{^k}_{_{qd^{^{-1}}}}(\uafa)$ for $k=1,2,...,n$ (see remark 2.2).
But from (3.5.1) these generators are equal respectively to
$\psi^{^k}\circ\eb_i(z)\circ\psi^{^{-k}},$
$\psi^{^k}\circ\fb_i(z)\circ\psi^{^{-k}}$ and
$\psi^{^k}\circ\kb_i^{^{\pm}}(z)\circ\psi^{^{-k}}$ ($i\in\set{n+1}$).
Thus we are done.
The integrability of $\modulo$ follows from the integrability of the
$\uf$-module $\V$. For the central charge see remarks 3.3 and 3.5.

\cvd

\vskip 5mm

{\bf 4. Definition of an inverse functor. }

\vskip 3mm

Fix $c,d,q\in\CC^\times$, $l,n\in\NN$, and set $y=1$.

\vskip 3mm

{\bf 4.1.} \rems
{\bf (i)} Suppose that $q$ is not a root of unity. Then
$\hf$-modules and integrable $\uf$-modules are direct sums of
finite dimensional modules (see [L1, 6.3.6] for the $\uf$-case). Thus,
if  $l\leq n$,
the Schur duality in the finite case (see [J]) gives indeed an equivalence
between the category of $\hf$-modules and the category of integrable
$\uf$-modules of level $l$.

{\bf (ii)} Similarly, if $q$ is not a root of unity and $l\leq n$,
the affine Schur duality gives indeed
an equivalence between the category of $\haf$-modules and the category
of integrable $\uaf$-module with trivial central charge and level $l$.

\vskip3mm

{\bf 4.2. Theorem.} {\sl Suppose that $x=d^{^{-n-1}}q^{^{n+1}}$,
$c=y=1$, $l+1<n$ and $q$ is not a root of unity.
Let $M'$ be an integrable $\utor$-module
with trivial central charge and level $l$. Then there exists a $\htor$-module,
$M$, such that $M'\simeq\modulo$ as $\utor$-modules.}

\vskip3mm

{\sl Proof.}  Given $i=1,2$, the restriction of $M'$ to
$\uaf^{^{(i)}}$ is integrable with trivial central charge
and level $l$. Thus, by affine Schur duality (see remark 4.1 (ii)) one gets an
$\dot\Hb^{^{(i)}}$-module, $M^{^{(i)}}$,
such that $M'\simeq M^{^{(i)}}\ot_{_{\hf}}{\bf V}^{^{\ot l}}$
as $\dot{\bf U}^{^{(i)}}$-modules.
Moreover, the $M^{^{(i)}}$ are isomorphic as $\hf$-modules. So just
denote them by $M$. By construction the action of $\eb_0,\fb_0,\kb_0$
is  as in remark 3.5
and the formulas in theorem 3.3 are valid,
where the action of $\X_j,\Y_j\in\htor$ on $M$ is given by the $\hafb$-module
and the $\hafa$-module structure of $M$. In order to prove that these
actions extend to a $\htor$-module structure it's sufficient to verify that
for any $m\in M$,
$$m\cdot\Q\cdot\Y_1\cdot\Q^{^{-1}}=m\cdot\Y_2\qquad{\rm and}\quad
m\cdot\Q\cdot\Y_l\cdot\Q^{^{-1}}=x\,m\cdot\Y_1
\leqno(4.2.1)$$

where $\Q=\X_1\cdot\T_{_{1,l-1}}$ (see 1.3).
But $M\mapsto\modulo$ is an
equivalence from the category of $\hf$-modules to the category
of integrable $\uf$-modules of level $l$
(since $q$ is not a root of unity, see 4.1).
Thus, if $\vb\in\V^{^{\ot l}}$ is a generator of $\V^{^{\ot l}}$, i.e.
$\V^{^{\ot l}}=\uf\cdot\vb$, the map $M\ni m\mapsto m\ot\vb\in\modulo$
is injective for any $\hf$-module $M$.

(i) Set $\vb=v_{_1}\ot v_{_2}\ot\cdots\ot v_{_l}$ and
$\wb=v_{_2}\ot v_{_3}\ot\cdots\ot v_{_l}\ot v_{_{n+1}}$. Then,
$$\eb_0(m\ot\vb)=q^{^{1-l}}m\cdot\Q\ot\wb.$$
Since $\wb$ is a generator of $\V$ the first relation in (4.2.1) follows
from $\eb_0\cdot\kb^{^\pm}_2(z)=\kb^{^\pm}_2(z)\cdot\eb_0$,
which gives
$$\eb_0(m\cdot\t_1^{^\pm}(q^{^n}d^{^2}z\Y_2)\ot\vb)=
\kb^{^\pm}_2(z)\cdot(m\cdot\X_1\ot v_{_{n+1}}\ot v_{_2}\ot\cdots\ot v_{_l}),$$
i.e.,
$$q^{^{1-l}}\,m\cdot\t_1^{^\pm}(q^{^n}d^{^2}z\Y_2)\cdot\Q\ot\wb=
q^{^{1-l}}\,m\cdot\Q\cdot\t_1^{^\pm}(q^{^n}d^{^2}z\Y_1)\ot\wb.$$

(ii) Set $\vb=v_{_{1}}\ot v_{_{3}}\ot v_{_{4}}\ot\cdots\ot v_{_{l+1}}$
and $\wb=v_{_3}\ot v_{_{4}}\ot\cdots\ot v_{_{l+1}}\ot v_{_{n+1}}$.
Suppose now that $n>l+1$. Relations (2.1.3) give
$$d^{^{-1}}(\eb_0\cdot\kb_{1,-1}-q^{^{-1}}\kb_{1,-1}\cdot\eb_0)\cdot\kb_1=
(q^{^{-1}}-q)\eb_{0,-1}=
d(\eb_0\cdot\kb_{n,-1}-q^{^{-1}}\kb_{n,-1}\cdot\eb_0)\cdot\kb_n.$$
When evaluated on $m\ot\vb$ it writes
$$d^{^{-1}}\eb_0\cdot\kb_{1,-1}\cdot\kb_1(m\ot\vb)=
-q^{^{-1}}d\,\kb_{n,-1}\cdot\eb_0\cdot\kb_n(m\ot\vb),$$
i.e.,
$$d^{^{n+1}}m\cdot\Q\cdot\Y_l\ot\wb=q^{^{n+1}}m\cdot\Y_1\cdot\Q\ot\wb.$$
Thus, since $\wb$ is a generator of $\V$, we get the
second relation in (4.2.1).
\cvd

\vskip1cm

{\bf 5. Conclusion.}

\vskip3mm

Using the theorems 3.5 and 4.2 we get the duality theorem :

\vskip3mm

\th {\sl Suppose that $q$ is not a root of unity, $x=d^{^{-n-1}}q^{^{
n+1}}$, $c=y=1$  and $l+1<n$. Then the functor
$M\mapsto\modulo$
is an equivalence between the category of right $\htor$-modules and
the category of left
integrable $\utor$-modules with trivial central charge and level $l$.}

\vskip3mm

\rem  By duality, since $\htor$ admits finite dimensional
representations if and only if $x=1$, the toroidal quantum group $\utor$
admits finite dimensional representations if and only if
$d^{^{n+1}}=q^{^{n+1}}$.

\vskip1cm

\centerline{\bf References}

\vskip3mm

\hangindent=0.7in[B]$\enspace$ J. Beck, Braid group action and quantum
affine algebras. {\sl Commun. Mathem. Phys.}, {\bf 165} (1994), 555-568.

\hangindent=0.7in[C1]$\enspace$ I. Cherednik, Double-affine Hecke algebras,
Knizhnik-Zamolodchikov equations, and MacDonald's operators.
{\sl Intern. Math. Research Notices, Duke Math. J.}, {\bf 68} (1992),
171-180.

\hangindent=0.7in[C2]$\enspace$ I. Cherednik, A new interpretation of
Gelfand-Tzetlin basis. {\sl Duke Math. J.}, {\bf 54} (1987), 563-578.

\hangindent=0.7in[CP]$\enspace$ V. Chari, A. Pressley, Quantum Affine Algebras
and Affine Hecke Algebras. {\sl Preprint q-alg9501003}.

\hangindent=0.7in [D1]$\enspace$ V. Drinfeld, A new realization of Yangians
and quantized affine algebras. {\sl Soviet. Math. Dokl.}, {\bf 36} (1988),
212-216.

\hangindent=0.7in[D2]$\enspace$ V. Drinfeld, Yangians and degenerate affine
Hecke algebras. {\sl Func. Anal. Appl.}, {\bf 20} (1986), 62-64.

\hangindent=0.7in[GKV]$\enspace$ V. Ginzburg, M. Kapranov, E. Vasserot,
Langlands reciprocity for algebraic surfaces.
{\sl Mathem. Research Letters}, {\bf 2} (1995), 147-160.

\hangindent=0.7in[GRV]$\enspace$ V. Ginzburg, N. Reshetikhin, E. Vasserot,
Quantum groups and flag varieties. {\sl Contemp. Mathem.}, {\bf 175} (1994),
101-130.

\hangindent=0.7in[J]$\enspace$ M. Jimbo, A $q$-analogue of $U(gl(N+1))$, Hecke
algebra and the Yang-Baxter equation. {\sl Lett. Math. Phys.}, {\bf 11} (1986),
247-252.

\hangindent=0.7in[L1]$\enspace$ G. Lusztig, Introduction to Quantum
groups. {\sl Progress in Mathem., Birkh\"auser}, {\bf 110} (1993).

\hangindent=0.7in[L2]$\enspace$ G. Lusztig, Affine Hecke algebras and
their graded version. {\sl Journal of the AMS}, {\bf 2} (1989), 599-625.

\bye